\input amstex

\documentstyle{amsppt}
\document
\topmatter
\title
Generalized orthotoric K\"ahler surfaces
\endtitle
\author
W\l odzimierz Jelonek
\endauthor
\abstract{ In this paper we describe   QCH K\"ahler surfaces $(M,g,J)$ of generalized  orthotoric type type. We introduce  a distinguished orthonormal frame on $(M,g)$  and give the structure equations for  $(M,g,J)$. In the case  when $I$  is conformally K\"ahler and  $(M,g,J)$ is not hyperk\"ahler we  integrate these  structure equations and rediscover the orthotoric   K\"ahler surfaces.  We also investigate the hyperk\"aler case.     We   prove  in a simple way that if  $(M,g,J)$  is a hyperk\"ahler surface with a degenerate   Weyl tensor  $W^-$  (i.e.  QCH  hyperk\"ahler surface)  then  among  all   hyperk\"aler structures  on   $(M,g)$   there exists  a K\"ahler structure  $J_0$  such   that   $(M,g,J_0)$  is of Calabi type or of orthotoric type.   }
 \endabstract

\endtopmatter
\define\G{\Gamma}

\define\DE{\Cal D^{\perp}}
\define\e{\epsilon}
\define\n{\nabla}
\define\om{\omega}

\define\r{\rightarrow}
\define\w{\wedge}
\define\k{\diamondsuit}
\define\th{\theta}
\define\p{\partial}
\define\a{\alpha}

\define\lb{\lambda}

\define\1{D_{\lb}}
\define\2{D_{\mu}}
\define\0{\Omega}

\define\De{\Cal D}

\define \E{\Cal E}
\bigskip
{\bf  1.  Introduction.} {\it  QCH K\"ahler surfaces} are K\"ahler surfaces  $(M,g,J)$
admitting a global, $2$-dimen\-sional, $J$-invariant distribution
$\De$ having the following property: The holomorphic curvature
$K(\pi)=R(X,J X,J X,X)$ of any $J$-invariant $2$-plane $\pi\subset
T_xM$, where $X\in \pi$ and $g(X,X)=1$, depends only on the point
$x$ and the number $|X_{\De}|=\sqrt{g(X_{\De},X_{\De})}$, where
$X_{\De}$ is the orthogonal projection of $X$ on $\De$  (see [J-1]).  Every  QCH K\"ahler surface admits an opposite almost Hermitian structure $I$  such that the Ricci tensor $\rho$ of $(M,g,J)$ is $I$-invariant and $(M,g,I)$ satisfies the second Gray condition: $$R(X,Y,Z,W)-R(IX,IY,Z,W)=R(IX,Y,IZ,W)+R(IX,Y,Z,IW).\tag {G2}$$In  [J-5], [J-6], [J-2], [J-3] the author has studied QCH semi-symmetric K\"ahler surfaces $(M,g,J)$ which admit an opposite Hermitian structure  and such that one of distributions $\De,\De^{\perp}$ is integrable. A general method of constructing  such manifolds was given there.   The aim of the present paper is to  investigate  the case where  none  of the distributions  $\De,\De^{\perp}$ is integrable.  In  [J-1]  we have proved  that   a  QCH  K\"ahler surface  whose  opposite  structure $I$  is Hermitian  is of  Calabi type or orthotoric type,  or  is  hyperk\"ahler.   In   the   present paper we find the distinguished orthonormal  frame  $\{E_1,E_2,E_3,E_4\}$  and the structure equations with respect  to this frame  for the  QCH  K\"ahler surfaces with  Hermitian opposite structure $I$  which are not of generalized Calabi type. We call such QCH surfaces the generalized orthotoric surfaces. In the case  when $I$  is conformally K\"ahler and  $(M,g,J)$ is not hyperk\"ahler we  integrate these  structure equations and rediscover the orthotoric   K\"ahler surfaces  (see [A-C-G], [J-1]).  We also investigate the hyperk\"aler case.     We   prove that if  $(M,g,J)$  is a hyperk\"ahler surface with a degenerate   Weyl tensor  $W^-$  (i.e.  QCH  hyperk\"ahler surface)  then  among  all   hyperk\"aler structures  on   $(M,g)$   there exists  a K\"ahler structure  $J_0$  such   that   $(M,g,J_0)$  is of Calabi type or of orthotoric type.  We also  describe the cases of hyperk\"ahler surfaces with only one Calabi structure,  with  only one  orthotoric structure,  with all structures orthotoric.  The   K\"ahler surfaces are studied by many authors  see  for example  [D-T].
\bigskip
{\bf  2.  Distinguished  frame on generalized orthotoric surfaces}   In this section we construct  special  orthonormal   frame for  generalized orthotoric surfaces and  evaluate the connection forms and structure equations. We  prove   that from the structure equations it follows that both almost Hermitian structures  $I,J$  are integrable. If  $\De$ is a $I$-invariant  distribution  on a Hermitian surface  $(M,g,I)$  then it defines  an almost  Hermitian structure  $J$,  which gives a different orientation  by   $J_{|\De}=I_{|\De},  J_{|\De^{\perp}}=-I_{|\De^{\perp}}$.   If   $(M,g,I)$  is a Hermitian  surface  then  the Lee form  $\th$  form   $(M,g,I)$  is defined by  $d\om_I=2\th\w\om_I$    where   $\om_i$  is the K\"ahler form of   $(M,g,I)$   i.e.   $\om_I(X,Y)=g(IX,Y)$. The next two Lemmas are proved in [J-5], [J-4].
\bigskip
{\bf Lemma A. }{\it Let $(M,g,I)$ be a Hermitian 4-manifold. Assume that  $|\n I| > 0$ on $M$. Then there exists a global oriented orthonormal basis $\{E_3,E_4\}$ of the nullity foliation $\De$, such that for any local orthonormal oriented basis $\{E_1,E_2\}$ of $\DE$
$$\n \0=\a(\th_1\otimes\Phi+\th_2\otimes\Psi), \tag 2.1$$
where  $\Phi=\th_1\w\th_3-\th_2\w\th_4,\ \Psi=\th_1\w\th_4+\th_2\w\th_3$, $\a=\pm\frac1{2\sqrt{2}}|\n J|$ and $\{\th_1,\th_2,\th_3,\th_4\}$ is a cobasis dual to $\{E_1,E_2,E_3,E_4\}$.  Moreover $\delta\0=-2\a\th_3$ and  the Lee form   of  $(M,g,I)$  is $\th=-\a\th_4$.}

A frame  $\{E_1,E_2,E_3,E_4\}$ we call a special frame for   $(M,g,I)$.

{\bf Lemma B.  }{\it Let $(M,g,I)$ be a Hermitian surface with $I$-invariant Ricci tensor. Let $\{E_1,E_2,E_3,E_4\}$ be a local orthonormal frame such that (2.1) holds. Then
$$\gather\G^3_{11}=\G^3_{22}= E_3\ln \a,  \tag a\\
\G^3_{44}=\G^4_{21}= -\G^4_{12}=- E_3\ln \a, \tag b\\
\G^3_{21}=-\G^3_{12},\ \G^4_{11}=\G^4_{22},  \tag c\\
-\G^3_{21}+\G^4_{22}=\a,\tag d\\
\G^4_{33}=-E_4\ln\a+\a, \tag e\endgather$$
where $\n_XE_i=\sum\om^j_i(X)E_j$ and $\G^i_{kj}=\om^i_j(E_k)$.}

{\bf Lemma  C.} {\it Let $\{ E_1'',E_2'',E_3',E_4'\}$ be a special frame on $(M,g,I)$   where   $I$  is a Hermitian structure.   Assume that  $\De$ is a $I$-invariant  distribution, such that    $\De$  is different   from  $\De_0=span\{E_1'',E_2''\},\De_1=span\{E_3',E_4'\}$.  Then there exists a special frame   $\{ E_1',E_2',E_3',E_4'\}$ and a function $\phi:M\r \Bbb R$ such that $$\De=span\{-\cos\phi E_3'+\sin\phi E_1', -\cos\phi E_4'+\sin\phi E_2'\},$$  $$\De^{\perp}=span\{\sin\phi E_3'+\cos\phi E_1',-\sin\phi E_4'-\cos\phi E_2'\}.$$    The K\"ahler structure $J$ is given  by   $$\gather  J(-\cos\phi E_4'+\sin\phi E_2')=-\cos\phi E_3'+\sin\phi E_1',\\ J(-\sin\phi E_4'-\cos\phi E_2')=\sin\phi E_3'+\cos\phi E_1'.\endgather$$}

\medskip
{\it Proof.}:  Let for  example  $E_3'+X\in \De$  where  $X\in \De_0$.  Then  $E_4'+IX\in \De$ and define a vector $E_1'$ by $\cos\phi E_3'+\sin\phi E_1'=\cos\phi(E_3'+X)$  where $\cos\phi=\frac1{||E_3'+X||}$.
  If $aE_3'+bE_4'+X\in\De$  then  $aE_4'-bE_3'+IX\in \De$ hence  $(a^2+b^2)E_3'+aX-bIX\in\De$ and we can apply the first case.  If $\De\ne \De_0, \De_1$  the one of these cases holds.  Note that  $X\ne 0$ since otherwise  $\De=\De_1$.

\bigskip

  Let  $E_1=-\cos\phi E_4'+\sin\phi E_2', E_3=-\cos\phi E_3'+\sin\phi E_1'$ and  $\De_+=\De=span(E_1,E_3)=ker(IJ-id)$. Then   $E_1=-\cos\phi E_4'+\sin\phi E_2', E_3=-\cos\phi E_3'+\sin\phi E_1'$,  and  $\De_-=\De^{\perp}=span(E_2,E_4)=ker(IJ+id)$ where  $$E_2=-\cos\phi E_2'-\sin\phi E_4', E_4=\cos\phi E_1'+\sin\phi E_3'$$ and   $$IE_1'=E_2', IE_3'=E_4', IE_1=-E_3,IE_2=E_4, JE_1=E_3,JE_2=E_4.$$

  We call   the frame  $\{E_1,E_2,E_3,E_4\}$   the distinguished frame for  a generalized orthotoric surface.

  \medskip
  {\it Remark.}  We give the geometric interpretation of the angle  $\phi$.   Note  that  the angle $\gamma$  between the  $I,J$  invariant distribution  $\Cal D=span\{E_1, E_3\}$ and $\Cal  D_0=span\{E_3',E_4'\}=span\{\th^\sharp,I\th^\sharp\}$  is  given  by  $\sin\gamma=|vol[E_1,E_3,E_3',E_4']|=\sin^2\phi$.   Hence
  $\phi=\arcsin(\sqrt{\sin\gamma})$.

  \medskip

We also have $$\gather E_1'=\cos\phi E_4+\sin\phi E_3, E_2'=-\cos\phi E_2+\sin\phi E_1,\\ E_3'=-\cos\phi E_3+\sin\phi E_4, E_4'=-\cos\phi E_1-\sin\phi E_2.\endgather$$
  Then $\th_1'=\cos\phi \th_4+\sin\phi \th_3, \th_2'=-\cos\phi \th_2+\sin\phi \th_1,\th_3'=\sin\phi \th_4-\cos\phi \th_3, \th_4'=-\sin\phi \th_2-\cos\phi \th_1$.  In particular
$\th_1'\w\th_2'+\th_3'\w\th_4'=-\th_1\w\th_3+\th_2\w\th_4,\th_1\w\th_3+\th_2\w\th_4=\cos(2\phi)(\th_1'\w\th_2'-\th_3'\w\th_4')+\sin(2\phi)(\th_1'\w\th_4'+\th_3'\w\th_2')$.

Let  $\om_1=\th_1\w\th_3,\om_2=\th_2\w\th_4$.  Then  $\om=\om_1+\om_2$ is a  K\"ahler form for  $(M,g,J)$ and $\overline{\om}=-\om_1+\om_2$  is a  K\"ahler form for $(M,g,I)$.

   Note that $d\om_2= a\th_2\w\om_1+b\th_1\w\om_2$.    Hence    $d\overline{\om}=2d\om_2=2(b\th_1-a\th_2)\w(\om_2-\om_1)$.   Consequently   the Lee form of  $(M,g,I)$ is $\th=-\alpha\th_4'=\alpha(\cos\phi\th_1+\sin\phi\th_2)$  and  $d\overline{\om}=2\th\w\overline{\om}=2\alpha\cos\phi\th_1\w\om_2-2\alpha\sin\phi\th_2\w\om_1$.  Hence if   $\De$ is not integrable  then   $\phi$ can take values in $\Bbb Z\frac{\pi}2$  in the set with an empty interior.

Since  $I$  is Hermitian  we have $\n_X\overline{\om}=IX\w\th+X\w I\th$,  $X=X_1E_1+X_2E_2+X_3E_3+X_4E_4, IX=-X_1E_3+X_2E_4+X_3E_1-X_4E_2$.  Thus
$$\gather\n_X\overline{\om}=
(-X_1\th_3+X_2\th_4+X_3\th_1-X_4\th_2)\w\alpha(\cos\phi\th_1+\sin\phi\th_2)\\+(X_1\th_1+X_2\th_2+X_3\th_3+X_4\th_4)\w \alpha(-\cos\phi\th_3+\sin\phi\th_4).\endgather$$

   Hence $$\gather\n_X\overline{\om}= \alpha(\cos\phi\th_4+\sin\phi\th_3)(\th_1\w\th_2+\th_3\w\th_4)+\\+\alpha(-\cos\phi\th_2+\sin\phi\th_1)(\th_1\w\th_4+\th_2\w\th_3).\endgather$$

  Since $$\n_X\overline{\om}= 2\om^4_1(X)(\th_1\w\th_2+\th_3\w\th_4)+2\om^1_2(\th_1\w\th_4+\th_2\w\th_3)$$ then   $$2\om^1_2=\alpha(-\cos\phi\th_2+\sin\phi\th_1)=2\om^3_4,  2\om^4_1= \alpha(\cos\phi\th_4+\sin\phi\th_3)=2\om^3_2,$$ where  $\n_{X}E_j=\om_j^i(X)E_i$.

$$\gather \n_{E_1'}E_1'=\sin\phi(\cos\phi E_3\phi E_3+\sin\phi\n_{E_3}E_3-\sin\phi E_3\phi E_4+\\\cos\phi\n_{E_3}E_4)+\cos\phi(\cos\phi E_4\phi E_3 +\sin\phi \n_{E_4}E_3-\sin\phi E_4\phi E_4+\cos\phi\n_{E_4}E_4) =\\ \sin\phi\cos\phi E_3\phi E_3+\sin^2\phi\n_{E_3}E_3-\sin^2\phi E_3\phi E_4+\sin\phi\cos\phi\n_{E_3}E_4\\+\cos^2\phi E_4\phi E_3+ \cos\phi\sin\phi \n_{E_4}E_3-\cos\phi\sin\phi E_4\phi E_4+\cos^2\phi\n_{E_4}E_4.\endgather$$

Since from  Lemma B  we have $\G^{3'}_{11}=-\cos\phi E_4\phi-\sin\phi E_3\phi=E_3'\ln\alpha$  then   $$E_4(\alpha\sin\phi)=E_3(\alpha\cos\phi).$$

Similarly

$$\gather \n_{E_2'}E_2'=\sin\phi(\cos\phi E_1\phi E_1+\sin\phi\n_{E_1}E_1+\sin\phi E_1\phi E_2\\-\cos\phi\n_{E_1}E_2)-\cos\phi(\cos\phi E_2\phi E_1 +\sin\phi \n_{E_2}E_1+\sin\phi E_2\phi E_2-\cos\phi\n_{E_2}E_2) =\\ \sin\phi\cos\phi E_1\phi E_1+\sin^2\phi\n_{E_1}E_1+\sin^2\phi E_1\phi E_2-\sin\phi\cos\phi\n_{E_1}E_2\\-\cos^2\phi E_2\phi E_1 - \cos\phi\sin\phi \n_{E_2}E_1-\cos\phi\sin\phi E_2\phi E_2+\cos^2\phi\n_{E_2}E_2.\endgather$$

Hence ($\G^k_{ij}=\om_j^k(E_i), \G^{k'}_{ij}=(\om_j^{k})'(E_i')$ )

$$\gather \G^{3'}_{22}=\sin\phi\cos\phi(\cos\phi\G^4_{22}-\sin\phi\G^4_{12})-\sin\phi\cos\phi(-\cos\phi\G^3_{21}\\+\sin\phi\G^3_{11})=\frac1\alpha(\sin\phi E_4\alpha-\cos\phi E_3\alpha).\endgather$$

  Thus $\G^{4'}_{22}=\frac12\alpha-\sin\phi E_1\phi+\cos\phi E_2\phi$   and   $\G^{4'}_{11}=\frac12\alpha-\sin\phi\cos\phi(\sin\phi\G^1_{33}+\cos\phi\G^1_{43})-\sin\phi\cos\phi(\sin\phi\G^2_{34}+\cos\phi\G^2_{44})$.  Consequently

$\sin\phi\cos\phi(\sin\phi\G^1_{33}+\cos\phi\G^1_{43})+\sin\phi\cos\phi(\sin\phi\G^2_{34}+\cos\phi\G^2_{44})=\sin\phi E_1\phi-\cos\phi E_2\phi$.

$$\gather \n_{E_1'}E_2'=\sin\phi(\cos\phi E_3\phi E_1+\sin\phi\n_{E_3}E_1+\sin\phi E_3\phi E_2-\\ \cos\phi\n_{E_3}E_2)+\cos\phi(\cos\phi E_4\phi E_1 +\sin\phi \n_{E_4}E_1+\sin\phi E_4\phi E_2-\cos\phi\n_{E_4}E_2) =\\ \sin\phi\cos\phi E_3\phi E_1+\sin^2\phi\n_{E_3}E_1+\sin^2\phi E_3\phi E_2-\sin\phi\cos\phi\n_{E_3}E_2\\+\cos^2\phi E_4\phi E_1+\sin\phi\cos\phi\n_{E_4}E_1+\cos\phi\sin\phi E_4\phi E_2-\cos^2\phi\n_{E_4}E_2.\endgather$$

$\G^{4'}_{12}=-\sin\phi E_3\phi-\cos\phi E_4\phi$.

$$\gather \n_{E_2'}E_1'=\sin\phi(\cos\phi E_1\phi E_3+\sin\phi\n_{E_1}E_3-\sin\phi E_1\phi E_4+\\ \cos\phi\n_{E_1}E_4)-\cos\phi(\cos\phi E_2\phi E_3 +\sin\phi \n_{E_2}E_3-\sin\phi E_2\phi E_4+\cos\phi\n_{E_2}E_4) =\\ \sin\phi\cos\phi E_1\phi E_3+\sin^2\phi\n_{E_1}E_3-\sin^2\phi E_1\phi E_4+\sin\phi\cos\phi\n_{E_1}E_4\\-\cos^2\phi E_2\phi E_3-\sin\phi\cos\phi\n_{E_2}E_3+\cos\phi\sin\phi E_2\phi E_4-\cos^2\phi\n_{E_2}E_4.\endgather$$

$\G^{4'}_{21}=-\cos\phi\sin\phi(\sin\phi\G^1_{13}-\cos\phi\G^1_{23})+\cos\phi\sin\phi(-\sin\phi\G^2_{14}+\cos\phi\G^2_{24})$

$\G^{3'}_{12}=\frac12\alpha-\cos\phi\sin\phi(\sin\phi\G^4_{32}+\cos\phi\G^4_{42})-\cos\phi\sin\phi(\sin\phi\G^3_{31}+\cos\phi\G^3_{41})$.

$\G^{3'}_{21}=-\frac12\alpha-\sin\phi E_1\phi+\cos\phi E_2\phi$.

$$\gather \n_{E_4'}E_4'=\cos\phi(\cos\phi E_1\phi E_2 +\sin\phi \n_{E_1}E_2-\sin\phi E_1\phi E_1+\\ \cos\phi\n_{E_1}E_1)+\sin\phi(\cos\phi E_2\phi E_2+\sin\phi\n_{E_2}E_2-\sin\phi E_2\phi E_1+ \cos\phi\n_{E_2}E_1) \\=-\sin\phi\cos\phi E_1\phi E_1+\sin^2\phi\n_{E_2}E_2-\sin^2\phi E_2\phi E_1+\sin\phi\cos\phi\n_{E_1}E_2\\+\cos^2\phi E_1\phi E_2+\sin\phi\cos\phi\n_{E_2}E_1+\cos\phi\sin\phi E_2\phi E_2+\cos^2\phi\n_{E_1}E_1.\endgather$$

$\G^{3'}_{44}=\sin^2\phi(\cos\phi\G^4_{12}+\sin\phi\G^4_{22})-\cos^2\phi(\cos\phi\G^3_{11}+\sin\phi\G^3_{21})$.

On the other hand we have  $\G^{3'}_{44}=-\G^{3'}_{22}$ and hence

$$\gather\sin\phi\cos\phi(\cos\phi\G^4_{22}-\sin\phi\G^4_{12})-\sin\phi\cos\phi(-\cos\phi\G^3_{21}+\sin\phi\G^3_{11})=\\
-\sin^2\phi(\cos\phi\G^4_{12}+\sin\phi\G^4_{22})+\cos^2(\cos\phi\G^3_{11}+\sin\phi\G^3_{21})\\=\frac1\alpha(\sin\phi E_4\alpha-\cos\phi E_3\alpha)\endgather$$   which implies  $$\sin\phi\G^4_{22}=\cos\phi\G^3_{11}.$$

$$\gather-\n_{E_3'}E_4'=\sin\phi(\cos\phi E_4\phi E_2+\sin\phi\n_{E_4}E_2-\sin\phi E_4\phi E_1+\\ \cos\phi\n_{E_4}E_1)-\cos\phi(\cos\phi E_3\phi E_2 +\sin\phi \n_{E_3}E_2-\sin\phi E_3\phi E_1+\cos\phi\n_{E_3}E_1) =\\ \sin\phi\cos\phi E_4\phi E_2+\sin^2\phi\n_{E_4}E_2-\sin^2\phi E_4\phi E_1-\sin\phi\cos\phi\n_{E_3}E_2\\-\cos^2\phi E_3\phi E_2+\sin\phi\cos\phi\n_{E_4}E_1+\cos\phi\sin\phi E_3\phi E_1-\cos^2\phi\n_{E_3}E_1.\endgather$$

$$-\G^{1'}_{34}=\cos\phi\sin\phi(\sin\phi\G^3_{41}-\cos\phi\G^3_{31})+\cos\phi\sin\phi(\sin\phi\G^4_{42}-\cos\phi\G^4_{32}).$$

$$\gather \n_{E_3'}E_1'=\sin\phi(\cos\phi E_4\phi E_3+\sin\phi\n_{E_4}E_3-\sin\phi E_4\phi E_4+\\\cos\phi\n_{E_4}E_4)-\cos\phi(\cos\phi E_3\phi E_3 +\sin\phi \n_{E_3}E_3-\sin\phi E_3\phi E_4+\cos\phi\n_{E_3}E_4) =\\ \sin\phi\cos\phi E_4\phi E_3+\sin^2\phi\n_{E_4}E_3-\sin^2\phi E_4\phi E_4+\sin\phi\cos\phi\n_{E_4}E_4\\-\cos^2\phi E_3\phi E_3- \cos\phi\sin\phi \n_{E_3}E_3+\cos\phi\sin\phi E_3\phi E_4-\cos^2\phi\n_{E_3}E_4.\endgather$$

$$\gather -\n_{E_4'}E_3'=\cos\phi(\cos\phi E_1\phi E_4 +\sin\phi \n_{E_1}E_4+\sin\phi E_1\phi E_3-\\ \cos\phi\n_{E_1}E_3)+\sin\phi(\cos\phi E_2\phi E_4+\sin\phi\n_{E_2}E_4+\sin\phi E_2\phi E_3-\cos\phi\n_{E_2}E_3) \\=\sin\phi\cos\phi E_1\phi E_3+\sin^2\phi\n_{E_2}E_4+\sin^2\phi E_2\phi E_3+\sin\phi\cos\phi\n_{E_1}E_4\\+\cos^2\phi E_1\phi E_4-\sin\phi\cos\phi\n_{E_2}E_3+\cos\phi\sin\phi E_2\phi E_4-\cos^2\phi\n_{E_1}E_3.\endgather$$

Thus
$$-\G^{1'}_{43}=\cos\phi E_1\phi+\sin\phi E_2\phi.$$

and

$$\gather \cos\phi E_1\phi+\sin\phi E_2\phi=\cos\phi\sin\phi(\sin\phi\G^3_{41}-\cos\phi\G^3_{31})\\+\cos\phi\sin\phi(\sin\phi\G^4_{42}-\cos\phi\G^4_{32}).\endgather$$

Similarly

$$-\G^{2'}_{43}=\sin\phi\cos\phi(\cos\phi\G^3_{11}+\sin\phi\G^3_{21})+\sin\phi\cos\phi(\cos\phi\G^4_{12}+\sin\phi\G^4_{22}),$$
and

$$-\G^{2'}_{34}=-\sin\phi E_4\phi+\cos\phi E_3\phi.$$

Consequently

$$\gather \sin\phi\cos\phi(\cos\phi\G^3_{11}+\sin\phi\G^3_{21})+\sin\phi\cos\phi(\cos\phi\G^4_{12}+\sin\phi\G^4_{22})=\\-\sin\phi E_4\phi+\cos\phi E_3\phi.\endgather$$

$$\gather \n_{E_2'}E_3'= \sin\phi\cos\phi E_1\phi E_4+\sin^2\phi\n_{E_1}E_4+\sin^2\phi E_1\phi E_3-\sin\phi\cos\phi\n_{E_1}E_3\\-\cos^2\phi E_2\phi E_4-\sin\phi\cos\phi\n_{E_2}E_4-\cos\phi\sin\phi E_2\phi E_3+\cos^2\phi\n_{E_2}E_3.\endgather$$

$$\gather\G^{4'}_{23}=-\cos^2\phi(\sin\phi \G^3_{11}-\cos\phi\G^3_{21})+\sin^2\phi(\sin\phi \G^4_{12}-\cos\phi\G^4_{22})\\=E_1'\ln\alpha=\frac1\alpha(\sin\phi E_3\alpha+\cos\phi E_4\alpha).\endgather$$

$$\gather -\cos^2\phi(\sin\phi \G^3_{11}-\cos\phi\G^3_{21})+\sin^2\phi(\sin\phi \G^4_{12}-\cos\phi\G^4_{22})+ \\ \sin\phi\cos\phi(\cos\phi\G^3_{11}+\sin\phi\G^3_{21})+\sin\phi\cos\phi(\cos\phi\G^4_{12}+\sin\phi\G^4_{22})=\\ \frac1\alpha(\sin\phi E_3\alpha+\cos\phi E_4\alpha
-\alpha\sin\phi E_4\phi+\alpha\cos\phi E_3\phi).\endgather$$
which gives
$$\gather \cos^3\phi\G^3_{21}+\sin^3\phi \G^4_{12}+  \sin^2\phi\cos\phi\G^3_{21}+\sin\phi\cos^2\phi\G^4_{12}=\\ \frac1{\alpha}(E_3(\alpha\sin\phi )+E_4(\alpha\cos\phi) )\endgather$$
hence
$$ \cos\phi\G^3_{21}+\sin\phi \G^4_{12}=\frac1{\alpha}( E_3(\alpha\sin\phi )+E_4(\alpha\cos\phi )).$$

We also   get

$$\gather -\cos^2\phi(\sin\phi \G^3_{11}-\cos\phi\G^3_{21})+\sin^2\phi(\sin\phi \G^4_{12}-\cos\phi\G^4_{22})- \\ \sin\phi\cos\phi(\cos\phi\G^3_{11}+\sin\phi\G^3_{21})-\sin\phi\cos\phi(\cos\phi\G^4_{12}+\sin\phi\G^4_{22})=\\\frac1\alpha( \sin\phi E_3\alpha+\cos\phi E_4\alpha
+\alpha\sin\phi E_4\phi-\alpha\cos\phi E_3\phi).\endgather$$

 $$\gather \G^{3'}_{22}=\sin\phi\cos\phi(\cos\phi\G^4_{22}-\sin\phi\G^4_{12})-\sin\phi\cos\phi(-\cos\phi\G^3_{21}+\sin\phi\G^3_{11})=\\ \frac1\alpha(\sin\phi E_4\alpha-\cos\phi E_3\alpha)\endgather$$

$$\gather\G^{4'}_{23}=-\cos^2\phi(\sin\phi \G^3_{11}-\cos\phi\G^3_{21})+\sin^2\phi(\sin\phi \G^4_{12}-\cos\phi\G^4_{22})\\=E_1'\ln\alpha=\frac1\alpha(\sin\phi E_3\alpha+\cos\phi E_4\alpha).\endgather$$

The above equations imply

$\cos\phi\G^4_{22}-\sin\phi\G^4_{12}=-\frac{E_3\alpha}{\alpha\sin\phi},-\cos\phi\G^3_{21}+\sin\phi\G^3_{11}=-\frac{E_4\alpha}{\alpha\cos\phi}$.

$$\sin\phi\G^4_{22}=\cos\phi\G^3_{11}.$$

$$-\G^{2'}_{43}=\sin\phi\cos\phi(\cos\phi\G^3_{11}+\sin\phi\G^3_{21})+\sin\phi\cos\phi(\cos\phi\G^4_{12}+\sin\phi\G^4_{22}),$$
and

$$-\G^{2'}_{34}=-\sin\phi E_4\phi+\cos\phi E_3\phi.$$

$$\gather \sin\phi\cos\phi(\cos\phi\G^3_{11}+\sin\phi\G^3_{21})+\sin\phi\cos\phi(\cos\phi\G^4_{12}+\sin\phi\G^4_{22})=\\-\sin\phi E_4\phi+\cos\phi E_3\phi.\endgather$$

$$\gather -\sin^2\phi(\cos\phi\G^4_{12}+\sin\phi\G^4_{22})+\cos^2\phi(\cos\phi\G^3_{11}+\sin\phi\G^3_{21})\\=\frac1\alpha(\sin\phi E_4\alpha-\cos\phi E_3\alpha).\endgather$$

Thus

$\cos\phi\G^3_{11}+\sin\phi\G^3_{21}=-\frac{E_4\phi}{\cos\phi},\cos\phi\G^4_{12}+\sin\phi\G^4_{22}=\frac{E_3\phi}{\sin\phi}$.

Consequently

$\G^4_{22}=-\frac{E_3(\alpha\cos\phi)}{\alpha\sin\phi}, \G^4_{12}=\frac{E_3(\alpha\sin\phi)}{\alpha\sin\phi}, \G^3_{21}=\frac{E_4(\alpha\cos\phi)}{\alpha\cos\phi},\G^3_{11}=-\frac{E_4(\alpha\sin\phi)}{\alpha\cos\phi}$.

$\sin\phi\cos\phi(\sin\phi\G^1_{33}+\cos\phi\G^1_{43})+\sin\phi\cos\phi(\sin\phi\G^2_{34}+\cos\phi\G^2_{44})=\sin\phi E_1\phi-\cos\phi E_2\phi$.
$$ \cos\phi\G^3_{21}+\sin\phi \G^4_{12}=\frac1{\alpha}( E_3(\alpha\sin\phi )+E_4(\alpha\cos\phi )).$$

$$\gather \cos\phi E_1\phi+\sin\phi E_2\phi=\cos\phi\sin\phi(-\sin\phi\G^1_{43}+\cos\phi\G^1_{33})\\+\cos\phi\sin\phi(-\sin\phi\G^2_{44}+\cos\phi\G^2_{34}).\endgather$$

$$\gather\n_{E_3'}E_3'=\sin\phi(\cos\phi E_4\phi E_4+\sin\phi\n_{E_4}E_4+\sin\phi E_4\phi E_3-\\ \cos\phi\n_{E_4}E_3)-\cos\phi(\cos\phi E_3\phi E_4 +\sin\phi \n_{E_3}E_4+\sin\phi E_3\phi E_3-\cos\phi\n_{E_3}E_3) =\\ \sin\phi\cos\phi E_4\phi E_4+\sin^2\phi\n_{E_4}E_4+\sin^2\phi E_4\phi E_3-\sin\phi\cos\phi\n_{E_3}E_4\\-\cos^2\phi E_3\phi E_4-\sin\phi\cos\phi\n_{E_4}E_3-\cos\phi\sin\phi E_3\phi E_3+\cos^2\phi\n_{E_3}E_3.\endgather$$

Thus  $$-\G^{4'}_{33}=\cos^2\phi(-\sin\phi\G^1_{43}+\cos\phi\G^1_{33})+\sin^2\phi(\sin\phi\G^2_{44}-\cos\phi\G^2_{34})=-\alpha+E'_4\ln\alpha$$

and $$ \gather \cos^2\phi x+\sin^2\phi y=-\alpha-\cos\phi\frac{E_1\alpha}\alpha- \sin\phi\frac{E_2\alpha}\alpha\\
\cos\phi \sin\phi x-\cos\phi\sin\phi y=\cos\phi E_1\phi+\sin\phi E_2\phi\endgather$$

where  $x=-\sin\phi\G^1_{43}+\cos\phi\G^1_{33}, y=\sin\phi\G^2_{44}-\cos\phi\G^2_{34}$.  It means that

$$x= -\alpha-\cos\phi E_1\ln (\alpha\cos\phi)-\sin\phi E_2\ln (\alpha\cos\phi)$$
$$y= -\alpha-\cos\phi E_1\ln (\alpha\sin\phi)-\sin\phi E_2\ln (\alpha\sin\phi)$$

$$\gather\n_{E_1'}E_3'=\sin\phi(\cos\phi E_3\phi E_4+\sin\phi\n_{E_3}E_4+\sin\phi E_3\phi E_3-\\ \cos\phi\n_{E_3}E_3)+\cos\phi(\cos\phi E_4\phi E_4 +\sin\phi \n_{E_4}E_4+\sin\phi E_4\phi E_3-\cos\phi\n_{E_4}E_3) =\\ \sin\phi\cos\phi E_3\phi E_4+\sin^2\phi\n_{E_3}E_4+\sin^2\phi E_3\phi E_3-\sin\phi\cos\phi\n_{E_3}E_3\\+\cos^2\phi E_4\phi E_4+\sin\phi\cos\phi\n_{E_4}E_4+\cos\phi\sin\phi E_4\phi E_3-\cos^2\phi\n_{E_4}E_3.\endgather$$

$\G^{4'}_{13}=\cos^2\phi(\sin\phi\G^1_{33}+\cos\phi\G^1_{43})-\sin^2\phi(\sin\phi\G^2_{34}+\cos\phi\G^2_{44})=\frac1\alpha(-\sin\phi E_1\alpha+\cos\phi E_2\alpha)$.

$\sin\phi\cos\phi(\sin\phi\G^1_{33}+\cos\phi\G^1_{43})+\sin\phi\cos\phi(\sin\phi\G^2_{34}+\cos\phi\G^2_{44})=\sin\phi E_1\phi-\cos\phi E_2\phi$.

We get

$$ \gather \cos^2\phi x-\sin^2\phi y=\frac1\alpha(-\sin\phi E_1\alpha+\cos\phi E_2\alpha)\\
\cos\phi \sin\phi x+\cos\phi\sin\phi y=\sin\phi E_1\phi-\cos\phi E_2\phi\endgather$$

where $x=\sin\phi\G^1_{33}+\cos\phi\G^1_{43}, y=\sin\phi\G^2_{34}+\cos\phi\G^2_{44}$.

Hence  $$x=-\sin\phi E_1\ln(\alpha\cos\phi)+\cos\phi E_2\ln(\alpha\cos\phi)$$, $$y=\sin\phi E_1\ln(\alpha\sin\phi)-\cos\phi E_2\ln(\alpha\sin\phi).$$

$$\gather \sin\phi\G^1_{33}+\cos\phi\G^1_{43}=-\sin\phi E_1\ln(\alpha\cos\phi)+\cos\phi E_2\ln(\alpha\cos\phi)\\
\cos\phi\G^1_{33}-\sin\phi\G^1_{43}=-\alpha-\cos\phi E_1\ln (\alpha\cos\phi)-\sin\phi E_2\ln (\alpha\cos\phi)\endgather$$

and

$$\G^1_{33}=-\alpha\cos\phi-E_1(\ln(\alpha\cos\phi)), \G^1_{43}=\alpha\sin\phi+E_2(\ln(\alpha\cos\phi)).$$
Similarly

$$\gather \sin\phi\G^2_{44}-\cos\phi\G^2_{34}=-\alpha-\cos\phi E_1\ln (\alpha\sin\phi)-\sin\phi E_2\ln (\alpha\sin\phi)\\
\cos\phi\G^2_{44}+\sin\phi\G^2_{34}=\sin\phi E_1\ln(\alpha\sin\phi)-\cos\phi E_2\ln(\alpha\sin\phi)\endgather$$
thus
$$\G^2_{44}=-\alpha\sin\phi-E_2(\ln(\alpha\sin\phi)), \G^2_{34}=\alpha\cos\phi+E_1(\ln(\alpha\sin\phi)).$$

  We also have

  $\om^3_1=E_4(\ln(\alpha\cos\phi))\th_2-\frac{E_4(\alpha\sin\phi)}{\alpha\cos\phi}\th_1 +(\alpha\cos\phi+E_1(\ln(\alpha\cos\phi)))\th_3-(\alpha\sin\phi+E_2(\ln(\alpha\cos\phi)))\th_4$  and

$\om^4_2=E_3(\ln(\alpha\sin\phi))\th_1-\frac{E_3(\alpha\cos\phi)}{\alpha\sin\phi}\th_2 -(\alpha\cos\phi+E_1(\ln(\alpha\sin\phi)))\th_3+(\alpha\sin\phi+E_2(\ln(\alpha\sin\phi)))\th_4$,  $d\th_i=\sum \om_i^k\w\th_k$.

   Since  $E_3(\a\cos\phi)=E_4(\a\sin\phi)$ and   $d\th_i=\sum \om_i^k\w\th_k$  it implies
\vskip1cm
{\bf  Theorem 2.1.}  {\it  If  $\{E_1,E_2,E_3,E_4\}$ is a distinguished  frame  on a generalized orthotoric  K\"ahler surface $(M,g,J)$  with opposite Hermitian  structure $I$   such that $IE_1=-E_3,IE_2=E_4, JE_1=E_3,JE_2=E_4$ and  $\{\th_1,\th_2,\th_3,\th_4\}$ is a
dual co-frame then the following structure equations  hold:}

 $ 2.1a\ d\th_1=-\frac12\alpha\sin\phi\th_1\w\th_2+E_4(\ln\alpha\cos\phi)\th_2\w\th_3-E_3(\ln\alpha\cos\phi)\th_1\w\th_3
+(\frac32\alpha\sin\phi+E_2(\ln(\alpha\cos\phi)))\th_3\w\th_4$.
\medskip
$2.1b\ d\th_2=\frac12\alpha\cos\phi\th_1\w\th_2+E_3(\ln\alpha\sin\phi)\th_1\w\th_4-E_4(\ln\alpha\sin\phi)\th_2\w\th_4
-(\frac32\alpha\cos\phi+E_1(\ln(\alpha\sin\phi)))\th_3\w\th_4$
\medskip
$2.1c\ d\th_3=\frac12\alpha\sin\phi\th_2\w\th_3+E_4(\ln\alpha\cos\phi)\th_1\w\th_2+\alpha\cos\phi\th_2\w\th_4
+(\frac32\alpha\sin\phi+E_2(\ln(\alpha\cos\phi))\th_4\w\th_1+(\alpha\cos\phi+E_1(\ln(\alpha\cos\phi)))\th_1\w\th_3$
\medskip

$2.1d\ d\th_4=\alpha\sin\phi\th_1\w\th_3-E_3(\ln\alpha\sin\phi)\th_1\w\th_2+\frac12\alpha\cos\phi\th_1\w\th_4+(\frac32\alpha\cos\phi+E_1(\ln(\alpha\sin\phi)))\th_3\w\th_2
+(\alpha\sin\phi+E_2(\ln(\alpha\sin\phi))\th_2\w\th_4$

{\it  Additionally  the connections forms are:}

$2\om^1_2=\alpha(-\cos\phi\th_2+\sin\phi\th_1)=2\om^3_4$,

$2\om^4_1= \alpha(\cos\phi\th_4+\sin\phi\th_3)=2\om^3_2$,

$\om^3_1=E_4\ln(\alpha\cos\phi)\th_2-E_3(\ln\alpha\cos\phi)\th_1 +(\alpha\cos\phi+E_1(\ln(\alpha\cos\phi)))\th_3-(\alpha\sin\phi+E_2(\ln(\alpha\cos\phi)))\th_4$

$\om^4_2=E_3(\ln(\alpha\sin\phi))\th_1-E_4(\ln\alpha\sin\phi)\th_2 -(\alpha\cos\phi+E_1(\ln(\alpha\sin\phi)))\th_3+(\alpha\sin\phi+E_2(\ln(\alpha\sin\phi)))\th_4$

{\it and for a Lee  form  $\th$  of $(M,g,I)$  we have}
\medskip
 $ d\th=(E_2(\alpha\cos\phi)-E_1(\alpha\sin\phi))(\th_2\w\th_1+\th_3\w\th_4)+(E_4(\alpha\sin\phi)+E_3(\alpha\cos\phi))(\th_4\w\th_2+\th_3\w\th_1)
+(E_3(\alpha\sin\phi)-E_4(\alpha\cos\phi))(\th_3\w\th_2+\th_1\w\th_4)$.

\medskip
Hence we have
\vskip1cm
{\bf Theorem  2.2}  {\it  The  Lie brackets of  the   fields  $E_1,E_2,E_3,E_4$   are as follows:}

$[E_1,E_2]=  \frac12\alpha\sin\phi E_1-\frac12\alpha\cos\phi E_2-E_4(\ln\alpha\cos\phi) E_3+E_3(\ln\alpha\sin\phi) E_4,$

$[E_1,E_3]=  E_3(\ln\alpha\cos\phi) E_1-(\alpha\cos\phi+E_1(\ln\alpha\cos\phi)) E_3-\alpha\sin\phi E_4,$

$[E_1,E_4]=  -E_3(\ln\alpha\sin\phi) E_2+(\frac32\alpha\sin\phi+E_2(\ln\alpha\cos\phi)) E_3-\frac12\alpha\cos\phi E_4,$

$[E_2,E_3]=  -E_4(\ln\alpha\cos\phi) E_1+(\frac32\alpha\cos\phi+E_1(\ln\alpha\sin\phi)) E_4-\frac12\alpha\sin\phi E_3,$

$[E_2,E_4]=  E_4(\ln\alpha\sin\phi) E_2-(\alpha\sin\phi+E_2(\ln\alpha\sin\phi)) E_4-\alpha\cos\phi E_3,$

$[E_3,E_4]=  -(\frac32\alpha\sin\phi +E_2(\ln\alpha\cos\phi ))E_1+(\frac32\alpha\cos\phi+E_1(\ln\alpha\sin\phi)) E_2.$

\medskip
Let $ f=\frac32\alpha\sin\phi+E_2\ln(\alpha\cos\phi)$, $g=\frac32\alpha\cos\phi+E_1(\ln(\alpha\sin\phi)), k=E_4(\alpha\sin\phi)+E_3(\alpha\cos\phi),l= E_3(\alpha\sin\phi)-E_4(\alpha\cos\phi)$, $h=E_2(\alpha\cos\phi)-E_1(\alpha\sin\phi)$.
The distribution $span\{E_3,E_4\}$ is integrable if and only if  $f=g=0$

 $\frac32\alpha\sin\phi+E_2(\ln(\alpha\cos\phi))=0$ i $\frac32\alpha\cos\phi+E_1(\ln(\alpha\sin\phi))=0$

 which is equivalent to

   $E_1(\alpha\sin\phi)=-\frac32\alpha^2\cos\phi\sin\phi$  i  $E_2(\alpha\cos\phi)=-\frac32\alpha^2\cos\phi\sin\phi$.

The distribution $span\{E_1,E_2\}$ is integrable if and only if

 $E_3(\alpha\sin\phi)=E_4(\alpha\cos\phi)=0$.

{\bf  Theorem 2.3.}  { \it Assume that on $M$ there is given a co-frame $\{ \th_1,\th_2,\th_3,\th_4\}$ with a dual  frame $\{ E_1,E_2,E_3,E_4\}$ such that the structure equations  2.1a-2.1d  are satisfied  for some smooth functions  $\alpha,\phi$. Let  $g= \th_1\otimes \th_1+\th_2\otimes \th_2+\th_3\otimes \th_3+\th_4\otimes \th_4$  and $\om_J=\th_1\w \th_3+\th_2\w \th_4, \om_I=-\th_1\w \th_3+\th_2\w \th_4 $.  Then  $(M,g,J)$ is a K\"ahler surface,  both structure  $I,J$ are Hermitian}

\medskip

{\it  Proof.}
 Let $\e=\pm 1$.  Then
$$\gather  d\th_1+i\e d\th_3=-\frac12\alpha\sin\phi\th_1\w\th_2+E_4(\ln\alpha\cos\phi)\th_2\w\th_3\\-E_3(\ln\alpha\cos\phi)\th_1\w\th_3
+(\frac32\alpha\sin\phi+E_2(\ln(\alpha\cos\phi)))\th_3\w\th_4\\ + \e i\frac12\alpha\sin\phi\th_2\w\th_3+i\e E_4(\ln\alpha\cos\phi)\th_1\w\th_2+\e i\alpha\cos\phi\th_2\w\th_4
\\+\e i(\frac32\alpha\sin\phi+E_2(\ln(\alpha\cos\phi))\th_4\w\th_1+\e i(\alpha\cos\phi+E_1(\ln(\alpha\cos\phi)))\th_1\w\th_3= \\ \frac12\alpha\sin\phi\th_2\w(\th_1+i\e\th_3) \\-iE_4(\ln\alpha\cos\phi)\th_2\w(\th_1+i\e \th_3)+(\frac32\alpha\sin\phi+iE_2(\ln(\alpha\cos\phi)))\th_4(\th_1+i\e\th_3)\\+\e\alpha\cos\phi \th_2\w(\th_2+i\th_4)+(\alpha\cos\phi+E_1(\ln(\alpha\cos\phi)))\th_1\w(\th_1+i\e\th_3)\\+\e iE_3(\ln\alpha\cos\phi)\th_1\w(\th_1+\e i\th_3),\endgather$$

$$\gather d\th_2+ id\th_4=\frac12\alpha\cos\phi\th_1\w\th_2+\\E_3(\ln\alpha\sin\phi)\th_1\w\th_4-E_4(\ln\alpha\sin\phi)\th_2\w\th_4
-(\frac32\alpha\cos\phi+E_1(\ln(\alpha\sin\phi)))\th_3\w\th_4+\\ i(\alpha\sin\phi\th_1\w\th_3-E_3(\ln\alpha\sin\phi)\th_1\w\th_2+\frac12\alpha\cos\phi\th_1\w\th_4\\+(\frac32\alpha\cos\phi+E_1(\ln(\alpha\sin\phi)))\th_3\w\th_2
+(\alpha\sin\phi+E_2(\ln(\alpha\sin\phi))\th_2\w\th_4)\\=\frac12\alpha\cos\phi\th_1\w(\th_2+i\th_4)-i E_3(\ln\alpha\sin\phi)\th_1\w(\th_2+i\th_4)\\+i(\frac32\alpha\cos\phi+E_1(\ln(\alpha\sin\phi)))\th_3\w(i \th_4+\th_2)+\e\alpha\sin\phi\th_1\w(\th_1+i\e\th_3)\\+(\alpha\sin\phi+E_2\ln(\alpha\sin\phi))\th_2\w(\th_2+i\th_4)\\+
iE_4(\ln\alpha\sin\phi)\th_2\w(\th_2+i\th_4).\endgather$$

It means that both structures   $I,J$ are Hermitian.  It is easy to check that  $d\om_J=0$.   Hence  $(M,g,J)$  is a K\"ahler surface.$\k$

\vskip1cm

{\bf  Theorem 2.4.}  { \it Assume that on $M$  there is given a coframe $\{ \th_1,\th_2,\th_3,\th_4\}$  such that the structure equations  2.1a-2.1d  are satisfied. Let  $g= \th_1\otimes \th_1+\th_2\otimes \th_2+\th_3\otimes \th_3+\th_4\otimes \th_4$  and $\om_J=\th_1\w \th_3+\th_2\w \th_4, \om_I=-\th_1\w \th_3+\th_2\w \th_4 $.  Then  $(M,g,J)$ is a K\"ahler surface,  both structure  $I,J$ are Hermitian and the Ricci form of our K\"ahler surface is   $\rho=d(d^I\ln\tan \phi))$,  hence is $J,I$-invariant.  Consequently  $(M,g,J)$  is a QCH  K\"ahler  surface of generalized orthotoric type. The Lee  form  $\th$  of $(M,g,I)$ satisfies
\medskip
 $ d\th=(E_2(\alpha\cos\phi)-E_1(\alpha\sin\phi))(\th_2\w\th_1+\th_3\w\th_4)+(E_4(\alpha\sin\phi)+E_3(\alpha\cos\phi))(\th_4\w\th_2+\th_3\w\th_1)
+(E_3(\alpha\sin\phi)-E_4(\alpha\cos\phi))(\th_3\w\th_2+\th_1\w\th_4)$.}

\medskip
{\it Proof.}  Note that the connection forms are uniquely determined by the structure equations. Moreover

$$\gather \om^3_1+\om^4_2=-E_3(\alpha\cos\phi)\th_1 +E_4(\ln(\alpha\cos\phi))\th_2 +(E_1(\ln(\alpha\cos\phi)))\th_3\\-(E_2(\ln(\alpha\cos\phi)))\th_4+E_3(\ln(\alpha\sin\phi))\th_1\\ -(E_1(\ln(\alpha\sin\phi)))\th_3+(E_2(\ln(\alpha\sin\phi)))\th_4-E_4(\alpha\sin\phi)\th_2 = \\ -d^I\ln(\alpha\cos\phi))+d^I\ln(\alpha\sin\phi) =d^I\ln|\tan \phi|.\endgather$$
It means that the Ricci form $\rho=d(\om^3_1+\om^4_2)=d(d^I\ln\tan \phi))$ of   $(M,g,J)$ has eigendistributions   $\De_+,\De_-$.$\k$

If   $\phi$  would be constant, then from the formula on the Ricci form we get  $\rho=0$  hence our surface is hyperk\"ahler.

\bigskip
{\bf  Lemma D.}  {\it Assume that  $E_3\alpha=E_4\alpha=E_3\phi=E_4\phi=0$  or equivalently  $$E_3(\alpha\cos\phi)=E_4(\alpha\cos\phi)=E_3(\alpha\sin\phi)=E_4(\alpha\sin\phi)=0.$$ Then   $f=g=0$  or   $\phi$  is constant.}

\medskip

{\it Proof.}  If   $E_3\alpha=E_4\alpha=E_3\phi=E_4\phi=0$  then   $$[E_3,E_4](\alpha\cos\phi)=[E_3,E_4](\alpha\sin\phi)=0.$$  Hence

$-fE_1(\alpha\cos\phi)+gE_2(\alpha\cos\phi)=0$  and  $-fE_1(\alpha\sin\phi)+gE_2(\alpha\sin\phi)=0$.   Thus if $d\th=0$  we get  $E_1(\alpha\sin\phi)=E_2(\alpha\cos\phi)$ and

$-fE_1(\alpha\cos\phi)+gE_1(\alpha\sin\phi)=0$,  $-fE_2(\alpha\cos\phi)+gE_2(\alpha\sin\phi)=0$

$E_1(\alpha)(-f\cos\phi+g\sin\phi)+\alpha E_1\phi(f\sin\phi+g\cos\phi)=0$,  $E_2(\alpha)(-f\cos\phi+g\sin\phi)+\alpha E_2\phi(f\sin\phi+g\cos\phi)=0$  which means that

$\alpha E_1\phi(f\sin\phi+g\cos\phi)=0$ and  $\alpha E_2\phi(f\sin\phi+g\cos\phi)=0$.   Thus  $f=g=0$  or   $E_1\phi=E_2\phi=0$.  In the last case  $(M,g,J)$ would be hyperk\"ahler  i.e. $\rho=0$.$\k$

\bigskip

{\bf   3. Generalized   orthotoric surfaces with  $d\th=0$.}

The case   $d\th=0$.
Let  $\th=-d\ln u$.  Then  $(M,g_1,I)$ is K\"ahler where  $g_1=u^{-2}g$.  Since  the Ricci tensor of  $(M,g,I)$  is $I$-invariant the field  $\xi=I\n^1u$ is Killing and $I$-holomorphic, hence the field $\n\frac1u$ is  holomorphic.
If   $W^+\ne 0$  then the field $\xi$ is also $J$-holomorphic.   Let us see what does it mean.
\vskip1cm
{\bf Theorem  3.1.}  {\it If   $(M,g,J)$ is a generalized orthotoric QCH surface with  $d\th=0$,  non-constant $\phi$ and a $J$-holomorphic Killing field $\xi$ then both distributions   $\De_+,\De_-$ are integrable and  $f=g=0$ and $E_3\alpha=E_4\alpha=E_3\phi=E_4\phi=0$.}

\medskip
{\it Proof.}   Note that $J\xi=\frac1uJI\theta^{\sharp}$ is a  $J$-holomorphic field.  Thus ($IE_1=-E_3,IE_2=E_4$)

$\n_XJ\xi=\frac1u(-\th(X)JI\theta^{\sharp}+\n_XJI\theta^{\sharp})$.  Hence
$J\n_{E_1}J\xi=\n_{E_3}J\xi$ i $\n_{E_3}J\xi=\frac1u(E_3(\a\cos\phi)E_1+\a\cos\phi\G^i_{31}E_i-E_3(\a\sin\phi)E_2-\a\sin\phi\G^i_{32}E_i)$ and

$$\gather \n_{E_1}J\xi=
\frac1u(-\a\cos\phi(JI\theta^{\sharp})\\+E_1(\a\cos\phi)E_1+\a\cos\phi\G^i_{11}E_i-E_1(\a\sin\phi)E_2-\a\sin\phi\G^i_{12}E_i).\endgather$$  If we look in the equation $J\n_{E_1}J\xi=\n_{E_3}J\xi$ on the second coordinate with respect to the frame $\{E_1,E_2,E_3,E_4\}$ we get

 $\a\cos\phi\G^2_{31}-E_3(\a\sin\phi)=-\a\cos\phi\G^4_{11}+\a\sin\phi\G^4_{12}$   thus

  $E_3(\a\sin\phi)=-E_3(\a\sin\phi)$ i $E_3(\a\sin\phi)=0$.  Similarly from  $J\n_{E_2}\xi=\n_{E_4}\xi$ we get $E_4(\a\cos\phi)=0$.

    It means that if $d\th=0$  then $f=g=0$ and both distributions are integrable. $\k$

Next we prove:
\vskip1cm
   {\bf  Theorem 3.2} {\it  If   $(M,g,J)$ is a generalized orthotoric QCH surface with  $d\th=0$, non-constant $\phi$   and a $J$-holomorphic Killing field $\xi$  then there exist local coordinates  $(x,y,z,t)$   such that   $E_1=h(x,y)a(x)\p_x, E_2=h(x,y)b(y)\p_y, E_3=\frac{hx}a\p_z+\frac{h}a\p_t, E_4=-\frac{hy}b\p_z-\frac{h}b\p_t  $ where $h=\frac1{\sqrt{x-y}}$ and   the metric $g$   is  $g= (x-y)(\frac1{a(x)^2}dx^2+\frac1{b(y)^2}dy^2)+\frac1{x-y}(a(x)^2(dz-ydt)^2+b(y)^2(dz-xdt)^2)$   for some smooth functions   $a,b:\Bbb R\rightarrow \Bbb R$.}

\medskip
   {\it  Proof.}  We have  $d\th_1\w\th_1=0=d\th_2\w\th_2$ hence  $\th_1=fdx,\th_2=gdy$  for some smooth functions  $f,g$. Take local  coordinates  $(z,t)$    such that   $ker dz\cap ker dt =span\{E_1,E_2\}$. Note that $f=f(x,y),g=g(x,y)$. It implies   $E_1=\frac1f\p_x,E_2=\frac1g\p_y$. We have $[E_1,E_2]=\frac1f\p_x\frac1g\p_y-\frac1g\p_y\frac1f\p_x =\frac12\a\sin\phi E_1-\frac12\a\cos\phi E_2$.   Consequently  $\p_y\ln f=\frac12\a\sin\phi g, \p_x\ln g=\frac12\a\cos\phi f$. Thus
  $$\p^2_{xy}\ln f=\p_x(\frac12\a\sin\phi g), \p^2_{xy}\ln g=\p_y(\frac12\a\cos\phi f).$$

   The condition $d\th=0$  implies $\p^2_{xy}\ln \frac fg=\p_x(\frac12\a\sin\phi g)-\p_y(\frac12\a\cos\phi f)=0$.  Hence    $\frac fg =b(y)\frac1{a(x)}$   and   $\frac 1f =ha,\frac 1g =hb(y)$,    where  $E_1=h(x,y)a(x)\p_x, E_2=h(x,y)b(y)\p_y$.

 The relation  $[E_1,E_2]=\frac12(\a\sin\phi E_1-\a\cos\phi E_2)=ha\p_x(hb)\p_y-hb\p_y(ha)\p_x=a\p_x E_2-b\p_yh E_1$  implies  $ \a\sin\phi=-2b\p_yh\ne0,\a\cos\phi=-2a\p_xh\ne0$

   Since $d\th=0$ we have $E_2(\a\cos\phi)=-\frac32\a^2\cos\phi\sin\phi$ and consequently $$-hb\p_y(2a\p_xh)=-6ab\p_xh\p_yh.$$  It follows  $\p^2_{xy}(\frac1{h^2})=-2\p_x(\frac{\p_yh}{h^3})=\frac{-2h\p^2_{xy}h+6\p_xh\p_yh}{h^4}=0$.     Hence  $\frac1{h^2}=F(x)+G(y)$   and changing coordinates we can assume  $h=\frac1{\sqrt{x-y}}$.

  We will show that there exist functions  $A(x,z,t),C(x,z,t),D(y,z,t),B(y,z,t)$   such that  $E_3=hA\p_z+hC\p_t, E_4=hB\p_z+hD\p_t$.

     For  example we have $E_4=B\p_z+D\p_t$   $[ah\p_x,B\p_z+D\p_t]=ah\p_xB\p_z+ah\p_xD\p_t=-\frac12\a\cos\phi(B\p_z+D\p_t)$. Thus  $ah\p_x\ln B=-\frac12\a\cos\phi=a\p_xh$ which implies $$\p_x\ln(\frac Bh)=0$$   and  $B=hB_1(y,z,t)$ etc.

 Since  $[E_2,E_4]=-(\a\sin\phi+E_2(\ln\a\sin\phi))E_4-\a\cos\phi E_3 $ we get  $[hb\p_y, hB\p_z+hD\p_t]=hb\p_y(hB)\p_z+hb\p_y(hD)\p_t$. On the other hand
 $[hb\p_y, hB\p_z+hD\p_t]=-(-2b\p_yh+bh\p_y\ln(b\p_yh))(hB\p_z+hD\p_t)+2a\p_xh(hA\p_z+hC\p_t)$. Hence
 $$b\p_y\ln(hB)=2b\p_y\ln h-b\p_y\ln(b\p_yh)+2a\p_x\ln h\frac AB$$ and

   $b\p_y\ln(B)=b\p_y\ln h-b\p_y\ln b-b\p_y\ln(\p_yh)+2a\p_x\ln h\frac AB$.
 Note that $\p_x\ln h=-\frac12\frac1{x-y}, \p_y\ln h=\frac12\frac1{x-y},\p_yh=\frac12(x-y)^{-\frac32},\p_y\ln |\p_yh|=-\frac3{2 (x-y)}, \p_xh=-\frac12(x-y)^{-\frac32},
\p_x \ln | \p_xh|=\frac3{2 (x-y)}$. Hence

 $b\p_y\ln(B)=\frac b{2(x-y)}-b\p_y\ln b-\frac32\frac b{2(x-y)}-\frac a{(x-y)}\frac AB$ and $Bb\p_y\ln(Bb)=-\frac{Bb+Aa}{x-y}$. The right side does not depend on $x$, hence
$Bb+Aa=c(y,z,t)(x-y)$. It follows $Aa=cx+d(z,t),Bb=-cy-d(z,t)$. Since $Aa$ does not depend on $y$ we get $c=c(z,t)$. Hence $A=\frac {cx+d}a,B=-\frac{cy+d}b$ if $c\ne 0$ and $A=\frac {l}a,B=-\frac{l}b$ where $l=l(z,t)$ if $c=0$.

The Killing field $\xi=\frac1{h^2}(-\a\cos\phi E_3+\a\sin\phi E_4)=\frac1{h^2}(2a\p_xhE_3-2b\p_yh E_4)=2a\p_x\ln h(A\p_z+C\p_t)-2b(\p_y\ln h(B\p_z+D\p_t)$ hence
$\xi=(2aA\p_x\ln h-2bB\p_y\ln h)\p_z+(2aC\p_x\ln h-2bD\p_y\ln h)\p_t= c_1(z,t)\p_z+c_2(z,t)\p_t$.

$d\th=\a\cos\phi\th_1+\a\sin\phi\th_2=-2a\p_xh\frac1{ah}dx-2b\p_yh\frac1{bh}dy=-2d\ln h=-d\ln u$ which implies $u=h^2$.

It follows $E_3=h\frac{c_1x+d_1}a\p_z+h\frac {c_2x+d_2}a\p_t,E_4=-h\frac {c_1y+d_1}b\p_z-h\frac {c_2y+d_2}b\p_t$  which means that  $E_3=\frac{hx}a\xi+\frac{hd_1}a\p_z+\frac{hd_2}a\p_t$
$E_4=-\frac{hy}b\xi-\frac{hd_1}b\p_z-\frac{hd_2}b\p_t$.  Let  $X=d_1\p_z+d_2\p_t$.   Since  $L_{\xi}\th_3=L_{\xi}\th_4=0$ we have $[\xi,X]=0$.    Similarly the Lie brackets of all fields   $\p_x,\p_y,\xi,X$ vanish.   It follows that we can find new coordinates $(x',y',z',t')$   such that  $\p_x=\p_{x'},\p_y=\p_{y'},\xi=\p_{z'},X=\p_{t'}$.  Note that we can choose $x=x',y=y'$.  In these new coordinates which we denote again by $(x,y,z,t)$  we have
 $E_3=\frac{hx}a\p_z+\frac{h}a\p_t, E_4=-\frac{hy}b\p_z-\frac{h}b\p_t$  and

 $$g= (x-y)(\frac1{a(x)^2}dx^2+\frac1{b(y)^2}dy^2)+\frac1{x-y}(a(x)^2(dz-ydt)^2+b(y)^2(dz-xdt)^2).$$

The   fields   $\xi=\p_z$ i $X=\p_t$ are holomorphic  Killing vector fields, their Lie bracket is $0$.$\k$

\medskip
{\it  Remark.}   Note that If   $(M,g,J)$ is a generalized orthotoric QCH surface with  $d\th=0$ which is not hyperK\"ahler  then $\phi$ is non-constant and $L_{\xi}\om_J$ is covariantly constant, hence   $ L_{\xi}\om_J=a\om_J$  for  $a\in\Bbb R$.  On the other hand  $g(\om_J,L_{\xi}\om_J)=0$  which means that  $\xi$  is $J$-holomorphic.
   In a hyperK\"ahler case  $d\th=0$ (if $W^-\ne0$ what we can assume since if $W^-=0$ then $(M,g,J)$ is $(\Bbb C^2,can)$) ) the field  $\xi=\frac1uI\theta^{\sharp}$ is $I$-holomorphic  Killing vector field, but $J\xi$ may not be $J$-holomorphic.  Integrability of  the distribution $\{E_1,E_2\}$  is a consequence of the fact that $\xi$ is   $J$- holomorphic and in general does not hold in this case.  But it is so if $\phi$ is constant.  The distribution span$\{E_3,E_4\}$   may not be integrable even if $\phi=const$.
\bigskip
{\bf  Theorem  3.3.}  {\it  Assume  that   $(M,g,J)$   is a generalized orthotoric surface with both integrable distributions $span\{E_1,E_2\},span\{E_3,E_4\}$.  Then $d\th=0$. }

\medskip
{\it Proof.}    We have

$E_3(\a\sin\phi)=E_4(\a\cos\phi)=0$,

 $E_4(\a\sin\phi)=E_3(\a\cos\phi)$,

 $E_1(\a\sin\phi)=-\frac32\a^2\sin\phi\cos\phi=E_2(\a\cos\phi)$.

 $[E_1,E_2]=  \frac12\alpha\sin\phi E_1-\frac12\alpha\cos\phi E_2,$

$[E_1,E_3]=  E_3(\ln\alpha\cos\phi) E_1-(\alpha\cos\phi+E_1(\ln\alpha\cos\phi)) E_3-\alpha\sin\phi E_4,$

$[E_1,E_4]=  -\frac12\alpha\cos\phi E_4,$

$[E_2,E_3]=  -\frac12\alpha\sin\phi E_3,$

$[E_2,E_4]=  E_4(\ln\alpha\sin\phi) E_2-(\alpha\sin\phi+E_2(\ln\alpha\sin\phi)) E_4-\alpha\cos\phi E_3,$

$[E_3,E_4]=  0.$

$d\th_1=-\frac12\alpha\sin\phi\th_1\w\th_2-E_3(\ln\alpha\cos\phi)\th_1\w\th_3$.
\medskip
$d\th_2=\frac12\alpha\cos\phi\th_1\w\th_2-E_4(\ln\alpha\sin\phi)\th_2\w\th_4$
\medskip
$d\th_3=\frac12\alpha\sin\phi\th_2\w\th_3+\alpha\cos\phi\th_2\w\th_4
+(\alpha\cos\phi+E_1(\ln(\alpha\cos\phi)))\th_1\w\th_3$
\medskip

$d\th_4=\alpha\sin\phi\th_1\w\th_3+\frac12\alpha\cos\phi\th_1\w\th_4
+(\alpha\sin\phi+E_2(\ln(\alpha\sin\phi))\th_2\w\th_4$

 Since $d\th_1\w\th_1=d\th_2\w\th_2=0$  then  $\th_1=fdx,\th_2=gdy$  for some functions $f,g$. It follows $E_1=\frac1f\p_x,E_2=\frac1g\p_y$.  Hence   $E_1=\frac1f\p_x,E_2=\frac1g\p_y$. Note that $[E_1,E_2]=\frac1f\p_x\frac1g\p_y-\frac1g\p_y\frac1f\p_x =\frac12\a\sin\phi E_1-\frac12\a\cos\phi E_2$.   Thus $\p_y\ln f=\frac12\a\sin\phi g, \p_x\ln g=\frac12\a\cos\phi f$. Hence
  $$\p^2_{xy}\ln f=\p_x(\frac12\a\sin\phi)g+\frac14\a^2\sin\phi\cos\phi fg,$$ $$\p^2_{xy}\ln g=\p_y(\frac12\a\cos\phi)f+\frac14\a^2\sin\phi\cos\phi fg.$$

  From the integrability of  $span\{E_1,E_2\}$ we get

  $\p^2_{xy}\ln \frac fg=\p_x(\frac12\a\sin\phi)g-\p_y(\frac12\a\cos\phi)f=fg(E_1(\a\sin\phi)-E_2(\a\cos\phi)=0$.  Hence    $\frac gf =a(x,z,t)\frac1{b(y,z,t)}$   and   $\frac 1f =ha(x,z,t),\frac 1g =hb(y,z,t)$,    where

    $E_1=h(x,y,z,t)a(x,z,t)\p_x$, $E_2=h(x,y,z,t)b(y,z,t)\p_y$.

    We will find a local coordinates
$(x,y,z,t)$ such that

$E_1=ha\p_x,
E_2=hb\p_y,
E_3=A\p_z+B\p_t,
E_4=C\p_z+D\p_t$.

Then  $\th=\frac{\a\cos\phi}{ha}dx+\frac{\a\sin\phi}{hb}dy$.   On the other hand $d\th=2E_3(\a\cos\phi)(\th_4\w\th_2+\th_3\w\th_1)$.  Hence:

$gE_2(\frac{\a\cos\phi}{ha})=fE_1(\frac{\a\sin\phi}{hb})$ and $E_4(\frac{\a\cos\phi}{ha})=E_3(\frac{\a\sin\phi}{hb})=0$.

It follows

$E_4(ha)=E_3(hb)=0$ i  $E_3\ln ha=-E_3\a\cos\phi, E_4(\ln hb)=-E_4(\a\sin\phi)$.

We also have

$E_3\ln\frac ab=-E_3(\a\cos\phi),  E_4\ln\frac ab=E_4(\a\sin\phi)$ and  $\a\sin\phi=-2\p_yhb,\a\cos\phi=-2\p_xha$.

$E_3(\ln\a\cos\phi)\th_3=\p_z\ln f dz+\p_t\ln f dt,
E_3(\ln\a\cos\phi)\th_4=\p_z\ln g dz+\p_t\ln g dt$.

$E_4\ln f=E_3\ln g=0$.   If $d\th=2k(\th_4\w\th_2+\th_3\w\th_1)$ then  $0=dk\w (\th_4\w\th_2+\th_3\w\th_1)$ i $k=const$.  Hence   $E_3(\ln\a\cos\phi)=\frac k{\a\cos\phi},  E_4(\ln \a\sin\phi)=\frac k{\a\sin\phi}$.

 $[E_1,E_2]=  \frac12\alpha\sin\phi E_1-\frac12\alpha\cos\phi E_2,$

$[E_1,E_3]=  \frac k{\a\cos\phi} E_1-(\alpha\cos\phi+E_1(\ln\alpha\cos\phi)) E_3-\alpha\sin\phi E_4,$

$[E_1,E_4]=  -\frac12\alpha\cos\phi E_4,$

$[E_2,E_3]=  -\frac12\alpha\sin\phi E_3,$

$[E_2,E_4]=  \frac k{\a\sin\phi} E_2-(\alpha\sin\phi+E_2(\ln\alpha\sin\phi)) E_4-\alpha\cos\phi E_3,$

$[E_3,E_4]=  0.$

It is not difficult to check that  such a system of equations is not integrable if $k\ne0$.$\k$

\vskip1cm
{\bf  Theorem 3.4.} {\it  There are no  QCH generalized orthotoric surfaces with  constant $\a$.}
\medskip
{\it Proof.}  Since  $C(E_1\otimes E_1+E_2\otimes E_2)\otimes(\th\otimes\th)=\a^2$ it follows
 $C(L_{E_3}E_1\otimes E_1+E_1\otimes L_{E_3}E_1+L_{E_3}E_2\otimes E_2+E_2\otimes L_{E_3}E_2)\otimes(\th\otimes\th)+2g(L_{E_3}\th,\th)=0$.

 Consequently $2E_3\ln(\a\cos\phi)E_1\otimes E_1-E_4\ln(\a\cos\phi)(E_2\otimes E_1+E_1\otimes E_2)(\th\otimes\th)+2g(k\th_1+l\th_2)=0$.

 Note that

 $E_3(\cos\phi)\cos\phi-E_4(\cos\phi)\sin\phi +(E_3(\cos\phi)+E_4(\sin\phi))\cos\phi+(E_3(\sin\phi)-E_4(\cos\phi))\sin\phi=0$.

 Let $a=E_4(\sin\phi)=E_3(\cos\phi)$.  Then  $E_3(\sin\phi)=-\cot \phi a, E_4(\cos\phi)=-\tan\phi a$. Thus

 $a\cos\phi+\sin\phi\tan\phi a+2a\cos\phi+(-\cot \phi a+\tan\phi a)\sin\phi=0$   which  gives   $\frac{2a}{\cos\phi}=0$  and $a=0$.   Hence  $E_3\phi=E_4\phi=0$.

 From   Lemma D it follows that  $f=g=0$   or   $\phi$   is constant.     In the first case  both distributions are integrable   and  $(M,g,J)$  is orthotoric thus $\alpha$   can not be constant.    The case of constant   $\alpha$, $\phi$   is not possible,   which  follows from  the Lie brackets of the fields  $E_1,E_2,E_3,E_4$.$\k$
\bigskip
{\bf  4. QCH  hyperk\"ahler surfaces}.  QCH   hyperk\"ahler surfaces are exactly hyper\-k\"ahler surfaces with degenerate  Weyl tensor  $W^-$.  Note that we have

$\n_X(\th_1\w\th_2-\th_3\w\th_4)=(\om_3^1+\om_4^2)(X)(\th_1\w\th_4+\th_3\w\th_2)=d^I\ln \tan\phi(X)(\th_1\w\th_4+\th_3\w\th_2)$,

$\n_X(\th_1\w\th_4+\th_3\w\th_2)=-(\om_3^1+\om_4^2)(X)(\th_1\w\th_2-\th_3\w\th_4)=-d^I\ln \tan\phi(X)(\th_1\w\th_2-\th_3\w\th_4)$.

We will find the conditions when   $\gamma=\cos\psi (\th_1\w\th_2-\th_3\w\th_4)+\sin\psi (\th_1\w\th_4+\th_3\w\th_2)$ is  K\"ahler.

Then  $X\cos\psi (\th_1\w\th_2-\th_3\w\th_4)+\cos\psi d^I\ln \tan\phi(X)(\th_1\w\th_4+\th_3\w\th_2)+X\sin\psi(\th_1\w\th_4+\th_3\w\th_2)-\sin\psi d^I\ln \tan\phi(X)(\th_1\w\th_2-\th_3\w\th_4)=0$  which gets   $X\cos\psi=-\sin\psi d^I\ln \tan\phi(X)$  and  $X\sin\psi=\cos\psi d^I\ln \tan\phi(X)$.  Thus

$d\psi= d^I\ln \tan\phi$.

A hyperK\"ahler surface is QCH if and only if $$M=(\Bbb C^2,can)$$  or $M$ admits a Hermitian  structure $I$ in the inverse orientation.

In  fact $M$ is Einstein, if admits a Hermitian  structure $I$ in the inverse orientation then the tensor $W^-$ is degenerate.   If $W^-=0$ then $M=(\Bbb C^2,can)$.
 In the other case $W^-$ is everywhere different from $0$ and has a double eigenvalue.  On the other hand if  $W^-\ne0$ and $W^-$ is degenerate then $M$ admits an integrable almost Hermitian structure whose K\"ahler  form is an eigenform of $W^-$. In fact  $(M,g)$ is Einstein and then from the result of Derdzi/nski [D] $(M,g)$ admits if $W^-\ne 0$ a Hermitian structure  $I$ in an opposite orientation with a Lee form  $\th$ satisfying  $d\th=0$ and such that  $\om_I$ is a simple eigenvalue of $W^-$.  Note that  $W^-\ne0$  on the whole of $M$  or $W^-=0$  everywhere.   In the last  case  $M$ is  $\Bbb C^2$  with the flat metric.

 In fact $(\Bbb C^2,can)$ admits  integrable  Hermitian structures in an opposite orientation (local).

If $(M,g,J)$  is a hyperk\"ahler surface with $W^-\ne0$  then   $(M,g)$  admits a Hermitian structure $I$, such that  $\om_I$ is an eigenform of  $W^-$ corresponding to a simple eigenvalue  and  $I$  is conformally K\"ahler. The Killing vector field   $\frac1uI\theta^{\sharp}$ is $I$-holomorphic.
\vskip1cm

 { \bf Theorem 4.1.}   {\it   Let  $(M,g,J)$  be a generalized orthotoric  QCH  hyperK\"ahler surface.   Then such surface admits a K/"ahler structure $J_0$ of Calabi type,  i.e.  $(M,g,J_0)$ is of Calabi type,  if and only if  $\phi$  is constant.  }

\medskip
  {\it  Proof.}
 Assume that  one of K\"ahler structures is of Calabi type.   Then  $\om'=\th_1'\w\th_2'-\th_3'\w\th_4'=\cos2\phi(\th_2\w\th_4+\th_1\w\th_3)+\sin2\phi(\th_2\w\th_3-\th_1\w\th_4)$ is K\"ahler. Hence
$0=\n_X\om'= -2\sin2\phi X\phi(\th_2\w\th_4+\th_1\w\th_3)+2\cos2\phi X\phi(\th_2\w\th_3-\th_1\w\th_4)+\sin2\phi d^I\ln \tan\phi(X)(\th_1\w\th_2-\th_3\w\th_4)$
and consequently  $\phi=const$.  On the other hand if  $\phi$  is constant then   $0=\n_X\om'$.
Note that then all structures  $\th_2\w\th_4+\th_1\w\th_3,\th_2\w\th_3-\th_1\w\th_4, \th_1\w\th_2-\th_3\w\th_4$ are K\"ahler.

If   $(M,g,J)$  is   orthotoric   with   $\phi=const$ then  $M$ is hyperk\"ahler and admits  K\"ahler structure of Calabi type.
Then  $g= (x-y)(\frac1{a^2}dx^2+\frac1{b^2}dy^2)+\frac1{x-y}(a^2(dz-ydt)^2+b^2(dz-xdt)^2)$, where $a,b$ are constant.$\k$
\vskip1cm

 Assume now that $(M,g,J)$  is a hyperk\"ahler orthotoric surface with  $\phi\ne const$.  Then   $\sin\phi=\frac b{\sqrt{a^2+b^2}}, \cos\phi = -\frac a{\sqrt{a^2+b^2}}$, $\a=\sqrt{a^2+b^2}h^3$  and $tan\phi=-\frac ba$.  Thus  $\ln|\tan\phi|=\ln b-\ln a$.  It means  $d^I\ln \tan\phi= hb'(y)\th_4+ha'(x)\th_3$  and  $d^I\ln \tan\phi=\frac{bb'}{x-y}(dz-xdt)+\frac{aa'}{x-y}(dz-ydt)$.   Let $F=a^2,G=b^2$.  Then $d^I\ln \tan\phi=\frac{G'}{2(x-y)}(dz-xdt)+\frac{F'}{2(x-y)}(dz-ydt)=\frac{G'+F'}{2(x-y)}dz+\frac{-G'x-F'y}{2(x-y)}dt$. Let  $\gamma=\cos\psi (\th_1\w\th_2-\th_3\w\th_4)+\sin\psi (\th_1\w\th_4+\th_3\w\th_2)$ be  K\"ahler.
Since $d\psi=d^I\ln \tan\phi$  it means that $\psi$ is a function depending only on $z,t$ and consequently  $G'+F'=2c(x-y),  -G'x-F'y=2d(x-y)$ for some   $d,c\in \Bbb R$. It implies

$G'=-2cy -2a, F'=2cx+2a,a\in \Bbb R$ and   $G(y)=-cy^2-2ay+b_1,F(x)=cx^2+2ax+b_2$.  It means  $d=-a$.  Thus

$d\psi=cdz-adt$ and $\psi= cz-at+\phi_0$ where  $\phi_0\in\Bbb R$.  The Killing  vector field $\xi$ related with conformally K\"ahler structure  $I$  is $\xi=\frac\p{\p z}$.   This field is triholomorphic    if and only if   $d\psi(\xi)=c=0$. It follows that for   $G(y)=-2ay+b_1,F(x)=2ax+b_2$ all K\"ahler structures from the sphere of K\"ahler structures are orthotoric.    If $c\ne0$  then only  the original one given by  $G(y)=-cy^2-2ay+b_1,F(x)=cx^2+2ax+b_2$ is orthotoric, the other are not orthotoric or Calabi type. If  $c=0$  then all the other structures are orthotoric.   If $c\ne0$  then only  this one is orthotoric the other are not of Calabi type or orthotoric.

In passing we also proved that an orthotoric structure is hyperK\"ahler if and only if  $G(y)=-cy^2-2ay+b_1,F(x)=cx^2+2ax+b_2$.(it follows also from [A-C-G]).

Let  $(M,g,J)$  be hyperk\"ahler  QCH  surface.  Then

$E_3(\a\sin\phi)=E_4(\a\cos\phi)$

and

$E_3(\a\cos\phi)=E_4(\a\sin\phi)=0$.  Let  $\th=-d\ln u$.

When  the field $\xi$ is $J$-holomorphic and when is triholomorphic?

We have  $\beta(x)=g(\xi,X)=-\frac{\alpha\cos\phi}u\th_3+\frac{\alpha\sin\phi}u\th_4$.  Hence

$d(-\frac{\alpha\cos\phi}u\th_3+\frac{\alpha\sin\phi}u\th_4)=-\frac1ud(\alpha\cos\phi)\w\th_3+\frac{du}{u^2}\w(\alpha\cos\phi \th_3)-\frac{\alpha\cos\phi}ud\th_3+
\frac1ud(\alpha\sin\phi)\w\th_4-\frac{du}{u^2}\w(\alpha\sin\phi \th_4)+\frac{\alpha\sin\phi}ud\th_4$

It implies   $d\beta= \frac{-2}u((E_1(\alpha\cos\phi)+\alpha^2\cos^2\phi-\frac12\alpha^2\sin^2\phi)\th_1\w\th_3+
\frac{2}u((E_2(\alpha\sin\phi)-\frac12\alpha^2\cos^2\phi+\alpha^2\sin^2\phi)\th_2\w\th_4 -\frac2u(E_2(\alpha\cos\phi)+\frac32\alpha^2\sin\phi\cos\phi)\th_2\w\th_3
+\frac2u(E_1(\alpha\sin\phi)+\frac32\alpha^2\sin\phi\cos\phi)\th_1\w\th_4+\frac2uE_4(\alpha\cos\phi)(\th_3\w\th_4-\th_1\w\th_2)$.

 $d\beta= \frac{2}u(-E_1(\alpha\cos\phi)+E_2(\alpha\sin\phi)-\frac32\alpha^2\cos^2\phi+\frac32\alpha^2\sin^2\phi)(\th_1\w\th_3+\th_2\w\th_4)
-\frac2u(E_2(\alpha\cos\phi)+\frac32\alpha^2\sin\phi\cos\phi)(\th_2\w\th_3-\th_1\w\th_4)
+\frac2uE_4(\alpha\cos\phi)(\th_3\w\th_4-\th_1\w\th_2)-\frac{2}u((E_2(\alpha\sin\phi)+E_1(\alpha\cos\phi)+\frac{\alpha^2}2)(\th_1\w\th_3-\th_2\w\th_4) $

Hence $\xi$ is  $J$-holomorphic if  $E_2(\alpha\cos\phi)=-\frac32\alpha^2\sin\phi\cos\phi=E_1(\alpha\sin\phi)$ and   $E_4(\alpha\cos\phi)=E_3(\a\sin\phi)=0$  i.e.  if
$(M,g,J)$  is orthotoric.

The first three forms  are self-dual,  the last one is anti-self-dual.  In an orthotoric case  we have

$d\beta= \frac{2}u(-E_1(\alpha\cos\phi)+E_2(\alpha\sin\phi)-\frac32\alpha^2\cos^2\phi+\frac32\alpha^2\sin^2\phi)(\th_1\w\th_3+\th_2\w\th_4)
 -\frac{2}u(E_2(\alpha\sin\phi)+E_1(\alpha\cos\phi)+\frac{\alpha^2}2)(\th_1\w\th_3-\th_2\w\th_4) $
which means that $\xi$ is   $I,J$-holomorphic.    $\xi$ is triholomorphic  if

$-E_1(\alpha\cos\phi)+E_2(\alpha\sin\phi)-\frac32\alpha^2\cos^2\phi+\frac32\alpha^2\sin^2\phi=0$.

In a general hyperk\"ahler generalized  orthotoric case $\xi$  is   triholomorphic if:

$-E_1(\alpha\cos\phi)+E_2(\alpha\sin\phi)-\frac32\alpha^2\cos^2\phi+\frac32\alpha^2\sin^2\phi=0$

$E_2(\alpha\cos\phi)=E_1(\alpha\sin\phi)=-\frac32\alpha^2\sin\phi\cos\phi$

$E_4(\alpha\cos\phi)=E_3(\alpha\sin\phi)=0$   thus   $(M,g,J)$  is orthotoric   with $G(y)=-2ay+b_1,F(x)=2ax+b_2$.
\bigskip

   {\bf Theorem 4.2}  {\it  If  $\xi$  is Killing  on a hyperk\"ahler surface $(M,g)$ then there exists a complex structure  $J$ on $M$  such that  $\xi$ is  $J$-holomorphic.}

\medskip
{ \it Proof.}   (see also  [L], p. 226)  Let   $V= span \{J_1,J_2,J_3\}$ be a bundle  of K\"ahler structures on $M$.  If  $\phi_t$ is a  one-parametric  group of  transformations of $\xi$  then
$(\phi_t)_*(J)\in \Gamma(V)$.    Hence   $S^2\ni J\r (\phi_t)_*(J)\in S^2$  is a local group of isometries of $S^2$ (a geodesic in $SO(3)$)  which extends to a one-parametric group in $SO(3)$.   Hence after the changing a base in $V$  $\phi_t(J)= \exp tA J$   where  $ A=\pmatrix
0&0&0\\\cr 0&0&a\\\cr 0&-a&0\endpmatrix$.  It follows that this group has a fixed point, $(\phi_t)_*(J)=J$  and $L_{\xi}J=0$, which means that $\xi$ is $J$-holomorphic.$\k$
Hence we  have:
\vskip1cm

 {\bf Theorem 4.3}  {\it If   $(M,g,J)$  is a hyperk\"ahler surface with a degenerate  Weyl tensor  $W^-$ then  among K\"ahler structures on  $M$  there exists  one of orthotoric or of Calabi type.}

 \medskip
 {\it  Proof.}  Let   $W^-\ne0$.   Then   there exists a Hermitian structure  $I$  in an opposite orientation    such that  $\om_I$  is a simple eigenform of  $W^-$.  The Lee form  $\th$  of $(M,g,I)$ is closed,  locally  $d\th=-d\ln u$ and  the field  $\xi=\frac1u I\th^{\sharp}$  is Killing.   From   Th.4.2  it follows that  $\xi$  is $J_0$ holomorphic   for some K\"ahler  structure  $J_0$.  If   $\phi$  is non-constant   it follows that   $(M,g,J_0)$  is orthotoric.    If  $\phi$ is constant  then  there exists  a K\"ahler structure   $J$   such that   $(M,g,J)$  is of Calabi type.

\medskip
Now we give simple proof  of the following  theorem generalizing well known  result  of  Gibbens and Ruback (see [G-R]).
\vskip1cm
{\bf Theorem 4.4.}  {\it If  a hyperk\"ahler surface admits a two-dimensional subalgebra of  $\frak iso(M)$  then it admits a triholomorphic Killing vector field.}

\medskip

{\it Proof.}  Let  $\om_1,\om_2,\om_3$  be an orthogonal basis of  K\"ahler forms on $M$.  Then  (  since  $L_{\xi}(\n_X\om_i)=\n_X(L_{\xi}\om_i)=0$ if $\n\om_i=0$) we have
$L_{\xi}\om_i=\sum a_{ij}\om_j$,   $L_{\eta}\om_i=\sum b_{ij}\om_j$
where $a_{ij},b_{ij}\in\Bbb R$.

Since   $\om_i\w\om_j=0$  for $i\ne j$ we get
$L_{\xi}\om_i\w\om_j+\om_i\w L_{\xi}\om_j=0$ and   $\sum a_{ip}\om_p\w\om_j+\sum a_{jp}\om_i\w\om_p=0$  which implies  $a_{ij}+a_{ji}=0$.  Similarly   $L_{\xi}\om_i\w\om_i=L_{\xi}vol=0$ i $a_{ii}=b_{ii}=0$.

 Hence   $A,B\in \frak{so}(3)$ where
 $A=[a_{ij}], B=[b_{ij}]$. Let    $\Phi:\frak{iso}(M)\r \frak{so}(3) $  be defined by $\Phi(\eta)=B=[b_{ij}]$, where
$L_{\eta}\om_i=\sum b_{ij}\om_j$. Note that

$L_{\eta}L_{\xi}\om_i=\sum a_{ij}L_{\eta}\om_j=\sum a_{ij}b_{jp}\om_p$ i $L_{\xi}L_{\eta}\om_i=\sum b_{ij}L_{\xi}\om_j=\sum b_{ij}a_{jp}\om_p$.    Hence  $\Phi([\xi,\eta)=[\Phi(\xi),\Phi(\eta)]$.  Thus  $\Phi$ is a  homomorphism of  Lie  algebras and $im\Phi$ is a subalgebra of $\frak{so}(3)$ hence it has a dimension $0,1$ or $3$ and $\ker\Phi=\frak{Tri}(M)$ where $\frak{Tri}(M)$ is a subalgebra of triholomorphic Killing vector fields.  Consequently if there exists a subalgebra $\frak{h}\subset dim\frak{iso}(M)$ such that $dim\frak{h} =2$ then there exists a  triholomorphic Killing vector field on $M$. This field is a linear combination of two fields from $\frak{h}$.   In fact if  $X,Y\in\frak{h}$  then for example  $\Phi(X)=a\Phi(Y)$ and the field $X-aY$ is triholomorphic.
It holds if there exist two fields $\xi,\eta\in \frak{iso}(M)$ such that  $[\xi,\eta]=0$. It is enough to take $\frak{h}=span(\xi,\eta)$. Note that   an algebra $\frak{Tri}(M)=ker\Phi$ is an ideal in $\frak{iso}(M)$.$\k$
\medskip
{\it Remark}  From  the proof of  Th.4.4   follows another proof of Th.4.2.    Let   $\om=\sum\alpha_i \om_i$.    Then  $L_{\xi}\om=\sum_{i,j} \alpha_ia_{ij}L_{\eta}\om_j$.  Hence   $L_{\xi}\om=0$   if   $\sum_i\alpha_ia_{ij}=0$  for $j=1,2,3$.   It means that  $\Phi(\xi)(\alpha)=0$ where $\alpha=[\alpha_1,\alpha_2,\alpha_3]\in\Bbb R^3$.  Note that  $ker\Phi(\xi)\ne\{0\}$.   If we take  $\alpha\in ker\Phi(\xi)$  and $||\alpha||=1$  then  $J=\sum\alpha_iJ_i$, where  $\om_i(X,Y)=g(J_iX,Y)$, is a K\"ahler structure such that  $L_{\xi}J=0$.
\medskip
{\bf Corollary.} { \it  Hyperkahler   QCH  surface  admits local  triholomorphic  Killing vector field.  It means that every hyperk\"ahler surface with degenerate  Weyl tensor $W^-$  admits locally triholomorphic  Killing vector  field.}

\medskip
 {\it  Proof.}  Let  $(M,g)$  be a hyperk\"ahler  QCH surface. If   $(M,g,J)$ is an orthotoric surface then  from Th.4.4 it admits  a triholomorphic Killing vector field.
The  case of  Calabi type  hyperk\"ahler  surfaces   we  solve  in the case of Eguchi-Hanson metric.    The other cases are similar.

  Let  $M=\Bbb R\times SO(3)$ and $g=h(t)^2dt^2+f^2(t)\th_1^2+g(t)^2(\th_2^2+\th_3^2)$  where  $\th_1,\th_2,\th_3$  are left-invariant 1-forms on $SO(3)$  and  $d\th_i=2\e_{ijk}\th_j\w \th_k$.
Define  two complex  structures on $(M,g)$

$J(hdt)=f\th_1,  J(\th_2)=\th_3$ i  $J_-(hdt)=-f\th_1,  J_-(\th_2)=\th_3$.

Let  $X\in\frak{iso}(S^2)$  where   $p:M\r S^2$  and  $\th_2^2+\th_3^2=p^*m$   where $m$ is the metric on $S^2$.  We will find  conditions  on $\eta=k\xi+ X^-$ to be Killing.
We have  $g=h(t)^2dt^2+f^2(t)\th_1^2+g(t)^2p^*m$ and  $L_Xm=0$.  $L_{\eta}(f^2\th_1)= f^2d(k)\odot\th_1+f^2(X\lrcorner d\th_1)\odot\th_1$.

Hence we get an equation  $dk+X\lrcorner\om=0$ where $k$ is a function on $S^2$.  But  $d(X\lrcorner\om)=L_X\om=0$ hence the equation has a solution  ($H^1(S^2)=0$).  Hence
a  field $\eta$   is Killing if  $dk=-X\lrcorner\om$.    Note that   $[\xi,\eta]=0$  since   $g(\xi, X^-)=0$,  $g(\xi,[\xi,X^-])=0$ and the field $[\xi,X^-]$ is vertical.

Since $\frak{iso}(S^2)$ is three-dimensional we get three fields  $\eta_1,\eta_2,\eta_3$.  Note that since   $[\xi,\eta_i]=0$ we get $\Phi(\eta_i)=\alpha_i\Phi(\xi)$ for some  $\alpha_i\in \Bbb R$.   Hence the image of the subalgebra generated by

$\frak{h}=span\{\xi,\eta_1,\eta_2,\eta_3\}$

is one-dimensional, dim im$\Phi(\frak{h})=1$. Thus we have  three triholomorphic  vector fields on $M$.  Note also that Th.4.4 implies  that  if  $(M,g)$ does not admit a triholomorphic then $dim\frak{iso}(M)=0,1$ or $\frak{iso}(M)$  is isomorphic to  $\frak{so}(3)$.  None of these cases hold for a Calabi type K\"ahler surface.$\k$

\bigskip

\centerline{\bf References.}

\cite{A-C-G} V. Apostolov, D.M.J. Calderbank, P. Gauduchon {\it
The geometry of we\-akly self-dual K\"ahler surfaces},  Compos. Math.
135, 279-322, (2003)
\medskip
[D] Andrzej Derdzi\'nski {\it  Self-dual K\"ahler manifolds and Einstein manifolds of dimension four}, Compositio Mathematica, 49,(1983),405-433.

\par
\medskip
\cite{D-T}   Maciej Dunajski and Paul Tod {\it  Four-dimensional metrics conformal to K\"ahler }
Mathematical Proceedings of the Cambridge Philosophical Society, Volume 148, (2010), 485-503
\medskip
[G-R] Gibbons, Ruback   {\it The hidden symmetries of multi-centre metrics} Comm.  Math.  Phys. 115  (1988),267-300
\par
\medskip
\cite{J-1} W. Jelonek {\it K\"ahler surfaces with quasi constant holomorphic curvature}, Glasgow Math. J. 58, (2016), 503-512.
\par
\medskip
\cite{J-2} W. Jelonek {\it Semi-symmetric  K\"ahler surfaces }, Colloq. Math.  148, (2017), 1-12.
\par
\medskip
\cite{J-3} W. Jelonek {\it Complex foliations and  K\"ahler QCH surfaces },  Colloq. Math.  156, (2019), 229-242.
\par
\medskip
\cite{J-4} W. Jelonek {\it Einstein Hermitian and anti-Hermitian 4-manifolds }, Ann.  Po\-lon. Math. 81, (2003), 7-24.
\medskip
\cite{J-5} W. Jelonek {\it  QCH  K\"ahler surfaces II},  Journal of Geometry and Physics, (2020) 103735,
\medskip
\cite{J-6} W. Jelonek {\it Generalized  Calabi type  K\"ahler surfaces},  arxiv, (2020),
\medskip
[L]  Claude Lebrun  {\it  Explicit self-dual metrics on  $\Bbb{CP}^2\sharp...\sharp\Bbb{CP}^2$},  J.  Diff.  Geom  34,(1991),223-253.

\par
\medskip
Author adress:

 Institute of Mathematics

 Cracow  University of Technology

 Warszawska 24

 31-155 Krak\'ow,  POLAND.

\enddocument